\documentclass[12pt]{amsart}
\usepackage{times}
\usepackage{amssymb}
\usepackage{amscd}
\def\Aff{\mathop{\mathrm {Aff}}\nolimits}
\def\aff{\mathop{\mathrm {aff}}\nolimits}
\def\SL{\mathop{\mathrm {SL}}\nolimits}
\def\Id{\mathop{\mathrm {Id}}\nolimits}
\def\Ad{\mathop{\mathrm {Ad}}\nolimits}
\def\ad{\mathop{\mathrm {ad}}\nolimits}
\def\Ind{\mathop{\mathrm {Ind}}\nolimits}
\def\Id{\mathop{\mathrm {Id}}\nolimits}
\def\Aut{\mathop{\mathrm {Aut}}\nolimits}

\def\sgn{\mathop{\mathrm {sgn}}\nolimits}

\newtheorem{theorem}{Theorem}[section]
\newtheorem{proposition}[theorem]{Proposition}
\newtheorem{lemma}[theorem]{Lemma}
\newtheorem{remark}[theorem]{Remark}

\begin{document}
\title{Quantum Half-Planes via Deformation Quantization}
\author{Do Ngoc Diep and Nguyen Viet Hai}
\address{Institute of Mathematics, National Center for Natural Sciences and Technology, P. O. Box 631, Bo Ho, 10.000, Hanoi, Vietnam}
\email{dndiep@thevinh.ncst.ac.vn}
\thanks{This work was supported in part by the Vietnam National  Foundation for Fundamental Science Research}
\maketitle
\begin{abstract}
We demonstrate the main idea of constructing irreducible unitary representations of Lie groups by using Fedosov deformation quantization in the concrete case of the group $\Aff(\mathbb  R)$ of affine transformations of the real straight line. By an exact computation of the star-product and the operator $\hat{\ell}_Z$, we show that the resulting representations exhausted all the irreducible representations of this groups.
\end{abstract}
\section{Introduction}
Quantization normally means a procedure associating to each classical mechanical system some quantum systems, namely in the Heisenberg model or Schr\"odinger one. More precisely, the usual formulation of a quantization procedure is a correspondence associating to each symplectic manifold $(M,\omega)$  a Hilbert space ${\mathcal H}$ of so called quantum states and to each classical observable (i.e. each complex-valued function) $f$ a quantum observable (i.e. a normal operator) $Q(f)$, in such a way that the following relations hold
\begin{equation}
Q(1)=\Id_{\mathcal H}
\end{equation}
\begin{equation}
[Q(f),Q(g)] = \frac{\hbar}{i} Q(\{f,g\}) \label{rel2}
\end{equation}
To attack this general problem there are some approaches, such as Feynman path integral quantization, pseudo-differential operator quantization, Weyl quantization, geometric quantization, etc. ... Following the geometric quantization procedure, at first one restricts himself to consider the set of observables to be quantized and secondly interpret the geometric quantization procedure operators, see e.g.  \cite{dndiep},
$Q(f) := f + \frac{\hbar}{i}\nabla_{\xi_f}$ 
as operators up to the second order approximation in power of $\hbar$, satisfying the relation (\ref{rel2}). From this point of view the so called {\it Fedosov deformation quantization} can be viewed as higher order approximates of operators satisfying the relation (\ref{rel2}). The last interpretation is the main idea behind deformation quantization. This deformation quantization essentially differs from the geometric quantization initiated by A. Kirillov, B. Kostant and J.-P.Souriau, see \cite{arnalcortet1}, \cite{kostant}. 

Many mathematicians attempted to construct quantum objects related with classical ones: First it was created the so called Podles quantum spheres. Interpreted the classical upper half-plane as the principal affine space of the special linear group $\SL_2(\mathbb R)$, one introduced the quantum upper half-plane as some C*-algebra generated by some generators and relations. We concern this upper half-plane from another point of view.

It is well-known that co-adjoint orbits are homogeneous symplectic manifolds with respect to the natural Kirillov form on orbits. A natural question is to associate to these orbits some quantum systems, which could be called quantum co-adjoint orbits. In most general context, some quantum co-adjoint orbits appeared in \cite{arnalcortet1} - \cite{arnalcortet2}. Still it is difficult to calculate exactly the $\star$-product and the corresponding representations in concrete cases. In this paper we demonstrate such an idea for a concrete case of the group $\Aff(\mathbb  R)$ of affine transformations of the real straight line. The main difficulty is the fact that in the concrete case, we should find out explicit formulae. This group has only two nontrivial 2-dimensional orbits which are the upper and lower half-planes. We shall use the same notion of star-product, introduced by M. Flato and A. Lichnerowicz, see \cite{arnalcortet1}. Our main result is the fact that by an exact computation we can find out explicit star-product
 formula and then by using the Fedosov deformation quantization, the full list of irreducible unitary representations of this group. These results show effectiveness of the Fedosov quantization, what are unknown up-to-date.

We introduce some notations in \S2, in particular, the canonical
coordinates are found in Proposition 2.1. The operators $\hat{\ell}_Z$ which define the representation of the Lie algebra $\aff(\mathbb  R)$ are found in \S3. By exponentiating we obtain the corresponding unitary representation of Lie group $\Aff_0(\mathbb  R)$ in Theorem 4.2 of \S4.
  
\section{Canonical coordinates on the upper half-planes}
Recall that the Lie algebra $\mathfrak g = \aff(\mathbb  R)$ of affine transformations of the real straight line is described as follows, see for example \cite{dndiep}: The Lie group $\Aff(\mathbb R)$ of affine transformations of type $$x \in \mathbb R \mapsto ax + b, \mbox{ for some parameters }a, b \in \mathbb R, a \ne 0.$$ It is well-known that this group $\Aff(\mathbb R)$ is a two dimensional Lie group which is isomorphic to the group of matrices
$$\Aff(\mathbb R) \cong \{\left (\begin{array}{cc} a & b \\ 0 & 1 \end{array} \right) \vert a,b \in \mathbb R , a \ne 0 \}.$$ We consider its connected component $$G= \Aff_0(\mathbb R)= \{\left (\begin{array}{cc} a & b \\ 0 & 1 \end{array} \right) \vert a,b \in \mathbb R, a > 0 \}$$ of identity element. Its Lie algebra is
$$\mathfrak g = \aff(\mathbb R) \cong  \{\left (\begin{array}{cc} \alpha & \beta \\ 0 & 0 \end{array} \right) \vert \alpha, \beta  \in \mathbb R \}$$  admits a basis of two generators $X, Y$ with the only nonzero Lie bracket $[X,Y] = Y$, i.e. 
$$\mathfrak g = \aff(\mathbb R) \cong \{ \alpha X + \beta Y \vert [X,Y] = Y, \alpha, \beta \in \mathbb R \}.$$
The co-adjoint action of $G$ on $\mathfrak g^*$ is given (see e.g. \cite{arnalcortet2}, \cite{kirillov1}) by $$\langle K(g)F, Z \rangle = \langle F, \Ad(g^{-1})Z \rangle, \forall F \in \mathfrak g^*, g \in G \mbox{ and } Z \in \mathfrak g.$$ Denote the co-adjoint orbit of $G$ in $\mathfrak g$, passing through $F$ by 
$$\Omega_F = K(G)F :=  \{K(g)F \vert F \in G \}.$$ Because the group $G = \Aff_0(\mathbb R)$ is exponential (see \cite{dndiep}), for $F \in \mathfrak g^* = \aff(\mathbb R)^*$, we have 
$$\Omega_F = \{ K(\exp(U)F | U \in \aff(\mathbb R) \}.$$
It is easy to see that
$$\langle K(\exp U)F, Z \rangle = \langle F, \exp(-\ad_U)Z \rangle.$$ It is easy therefore to see that
$$K(\exp U)F = \langle F, \exp(-\ad_U)X\rangle X^*+\langle F, \exp(-\ad_U)Y\rangle Y^*.$$
For a general element $U = \alpha X + \beta Y \in \mathfrak g$, we have
$$\exp(-\ad_U) = \sum_{n=0}^\infty \frac{1}{n!}\left(\begin{array}{cc}0 & 0 \\ \beta & -\alpha \end{array}\right)^n = \left( \begin{array}{cc} 1 & 0 \\ L & e^{-\alpha} \end{array} \right),$$ where $L = \alpha + \beta + \frac{\alpha}{\beta}(1-e^\beta)$. This means that
$$K(\exp U)F = (\lambda + \mu L) X^* + (\mu e^{\-\alpha})Y^*. $$ From this formula one deduces  \cite{dndiep} the following description of all co-adjoint orbits of $G$ in $\mathfrak g^*$:
\begin{itemize}
\item If $\mu = 0$, each point $(x=\lambda , y =0)$ on the abscissa ordinate corresponds to a 0-dimensional co-adjoint orbit $$\Omega_\lambda = \{\lambda X^* \}, \quad \lambda \in \mathbb R .$$
\item For $\mu \ne 0$, there are two 2-dimensional co-adjoint orbits: the upper half-plane $\{(\lambda , \mu) \quad\vert\quad \lambda ,\mu\in \mathbb R , \mu > 0 \}$ corresponds to the co-adjoint orbit
\begin{equation} \Omega_{+} := \{ F = (\lambda + \mu L)X^* + (\mu e^{-\alpha})Y^* \quad \vert \quad \mu > 0 \}, \end{equation}
and the lower half-plane $\{(\lambda , \mu) \quad\vert\quad \lambda ,\mu\in \mathbb R , \mu < 0\}$ corresponds to the co-adjoint orbit
\begin{equation} \Omega_{-} := \{ F = (\lambda + \mu L)X^* + (\mu e^{-\alpha})Y^* \quad \vert \quad \mu < 0 \}. \end{equation}
\end{itemize}
We shall work from now on for the fixed co-adjoint orbit $\Omega_+$. The case of the co-adjoint orbit $\Omega_-$ is similarly treated. First we study the geometry of this orbit and introduce some canonical coordinates in it.
It is well-known from the orbit method \cite{kirillov1} that the Lie algebra $\mathfrak g = \aff(\mathbb R)$, realized by the complete right-invariant Hamiltonian vector fields on co-adjoint orbits $\Omega_F \cong G_F \setminus G$ with flat (co-adjoint) action of the Lie group $G = \Aff_0(\mathbb R)$. On the orbit $\Omega_+$ we choose a fix point $F=Y^*$. It is well-known from the orbit method that we can choose an arbitrary point $F$ on $\Omega_F$. It is easy to see that the stabilizer of this (and therefore of any) point  is trivial $G_F = \{e\}$. We identify therefore $G$ with $G_{Y^*}\setminus G$. There is a natural diffeomorphism $\Id_{\mathbb R} \times \exp(.)$ from the standard symplectic space $\mathbb R^2$ with symplectic 2-form $dp \wedge dq$ in canonical Darboux $(p,q)$-coordinates, onto the upper half-plane $\mathbb H_+ \cong \mathbb R \rtimes \mathbb R_+$ with coordinates $(p, e^q)$, which is, from the above coordinate description, also diffeomorphic to the co-adjoint orbit $\Omega_+$. We can use therefore $(p,q)$ as the standard canonical Darboux coordinates in $\Omega_{Y^*}$. There are also non-canonical Darboux coordinates $(x,y) = (p,e^q)$ on $\Omega_{Y^*}$. We show now that in these coordinates $(x,y)$, the Kirillov form looks like $\omega_{Y^*}(x,y) = \frac{1}{y}dx \wedge dy$, but in the canonical Darboux coordinates $(p,q)$, the Kirillov form is just the standard symplectic form $dp \wedge dq$. This means that there are  symplectomorphisms between the standard symplectic space $\mathbb R^2, dp \wedge dq)$, the upper half-plane $(\mathbb H_+, \frac{1}{y}dx \wedge dy)$ and the co-adjoint orbit $(\Omega_{Y^*},\omega_{Y^*})$. 
Each element $Z\in \mathfrak g$ can be considered as a linear functional $\tilde{Z}$ on co-adjoint orbits, as subsets of $\mathfrak g^*$, $\tilde{Z}(F) :=\langle F,Z\rangle$.  It is well-known that this linear function is just the Hamiltonian function associated with the Hamiltonian vector field $\xi_Z$, which represents $Z\in \mathfrak g$ following the formula 
$$(\xi_Zf)(x) := \frac{d}{dt}f(x\exp (tZ))|_{t=0}, \forall f \in C^\infty(\Omega_+).$$ 
The Kirillov form $\omega_F$ is defined by the formula 
\begin{equation}\label{7} \omega_F(\xi_Z,\xi_T) = \langle F,[Z,T]\rangle, \forall Z,T \in \mathfrak g = \aff(\mathbb R). \end{equation} This form defines the symplectic structure and the Poisson brackets on the co-adjoint orbit $\Omega_+$. For the derivative along the direction $\xi_Z$ and the Poisson bracket we have relation $\xi_Z(f) = \{\tilde{Z},f\}, \forall f \in C^\infty(\Omega_+)$. It is well-known in differential geometry that the correspondence
$Z \mapsto \xi_Z, Z \in \mathfrak g$ defines a representation of our Lie algebra by vector fields on co-adjoint orbits. If the action of $G$ on $\Omega_+$ is flat \cite{dndiep}, we have the second Lie algebra homomorphism from  strictly Hamiltonian right-invariant vector fields into the Lie algebra of smooth functions on the orbit with respect to the associated Poisson brackets.

Denote by $\psi$ the indicated symplectomorphism from $\mathbb R^2$ onto $\Omega_+$
$$(p,q) \in \mathbb R^2 \mapsto \psi(p,q):= (p,e^q) \in \Omega_+$$
\begin{proposition}
1. Hamiltonian function $f_Z = \tilde{Z}$ in canonical coordinates $(p,q)$ of the orbit $\Omega_+$ is of the form $$\tilde{Z}\circ\psi(p,q) = \alpha p + \beta e^q, \mbox{ if  } Z = \left(\begin{array}{cc} \alpha & \beta \\ 0 & 0 \end{array} \right).$$

2. In the canonical coordinates $(p,q)$ of the orbit $\Omega_+$, the Kirillov form $\omega_{Y^*}$ is just the standard form $\omega = dp \wedge dq$.
\end{proposition}
{\it Proof}. 1. Each element $F \in (\aff(\mathbb R))^*$ is of the form $F = xX^* + yY^*$. This means that the value of the function $f_Z = \tilde{Z}$ on the element $Z = \alpha X + \beta Y$ is
$$\tilde{Z}(F) = \langle F,Z\rangle = \langle xX^* + yY^*, \alpha X + \beta Y \rangle = \alpha x + \beta y.$$ It follows therefore that 
\begin{equation}\label{6}\tilde{Z}\circ\psi(p,q) = \alpha p + \beta e^q, \mbox{ if  } Z = \left(\begin{array}{cc} \alpha & \beta \\ 0 & 0 \end{array} \right).\end{equation}

2. In canonical Darboux coordinates$(p,q)$,  $F = pX^* + e^qY^* \in \Omega_+$, and for $Z = \left(\begin{array}{cc} \alpha_1 & \beta_1 \\ 0 & 0 \end{array}\right), T = \left(\begin{array}{cc} \alpha_2 & \beta_2 \\  0 & 0 \end{array}\right),$ we have
$$\langle F, [Z,T]\rangle = \langle pX^* + e^q Y^*, (\alpha_1\beta_2 - \alpha_2\beta_1)Y \rangle = (\alpha_1\beta_2 - \alpha_2\beta_1)e^q,$$  i.e. 
\begin{equation}\label{7}\omega_F(\xi_Z,\xi_T) = (\alpha_1\beta_2 - \alpha_2\beta_1)e^q. \end{equation}
Let us consider two vector fields 
$$\xi_Z = \alpha_1 \frac{\partial}{\partial q} - \beta_1 e^q\frac{\partial}{\partial p},$$ i.e. 
$$\xi_Z(f) = \{ \alpha_1 p + \beta_1 e^q, f\} = \alpha_1 \frac{\partial f}{\partial q} - \beta_1 e^q\frac{\partial f}{\partial p} $$  and $$\xi_T = \alpha_2 \frac{\partial}{\partial q} - \beta_2 e^q\frac{\partial}{\partial p},$$ i.e. 
$$\xi_Z(f) = \{ \alpha_2 p + \beta_2 e^q, f\} = \alpha_2 \frac{\partial f}{\partial q} - \beta_2 e^q\frac{\partial f}{\partial p}. $$  We have 
\begin{equation}\label{8}\xi_Z \otimes \xi_T = \alpha_1 \alpha_2 \frac{\partial}{\partial q}\otimes \frac {\partial}{\partial q} + (\alpha_1\beta_2 - \alpha_2\beta_1)e^q \frac{\partial}{\partial p}\otimes \frac {\partial}{\partial q} + \beta_1\beta_2 e^{2q} \frac{\partial}{\partial p}\otimes \frac {\partial}{\partial p} \end{equation}
From (\ref{7}) and (\ref{8}) we conclude that in the canonical coordinates the Kirillov form is just the standard symplectic form $\omega = dp \wedge dq$. 

\section{Computation of generators $\hat{\ell}_Z$}
Let us denote by $\Lambda$ the 2-tensor associated with the Kirillov standard form $\omega = dp \wedge dq$ in canonical Darboux coordinates. We use also the multi-index notation. Let us consider the well-known Moyal $\star$-product of two smooth functions $u,v \in C^\infty(\mathbb R^2)$, defined by
$$u \star v = u.v + \sum_{r \geq 1} \frac{1}{r!}(\frac{1}{2i})^r P^r(u,v),$$ where
$$P^r(u,v) := \Lambda^{i_1j_1}\Lambda^{i_2j_2}\dots \Lambda^{i_rj_r}\partial_{i_1i_2\dots i_r} u \partial_{j_1j_2\dots j_r}v,$$ with $$\partial_{i_1i_2\dots i_r} := \frac{\partial^r}{\partial x^{i_1}\dots \partial x^{i_r}}, x:= (p,q) = (p_1,\dots,p_n,q^1,\dots,q^n)$$ as multi-index notation. It is well-known that this series converges in the Schwartz distribution spaces $\mathcal S (\mathbb R^n)$. We apply this to the special case $n=1$. In our case we have only $x = (x^1,x^2) = (p,q)$. 
\begin{proposition}\label{3.1}
In the above mentioned canonical Darboux coordinates $(p,q)$ on the orbit $\Omega_+$, the Moyl $\star$-product satisfies the relation
$$i\tilde{Z} \star i\tilde{T} - i\tilde{T} \star i\tilde{Z} = i\widetilde{[Z,T]}, \forall Z, T \in \aff(\mathbb R).$$
\end{proposition}
{\it Proof}.
Consider the elements $Z = \alpha_1 X + \beta_1 Y$ and $T = \alpha_2 X + \beta_2 Y$. Then as said above the corresponding Hamiltonian functions are
$\tilde{Z} = \alpha_1 p + \beta_1 e^q$ and $\tilde{T} = \alpha_2 p + \beta_2 e^q$. It is easy then to see that $$P^0(\tilde{Z},\tilde{T}) = \tilde{Z}.\tilde{T},$$
$$P^1(\tilde{Z},\tilde{T}) = \{ \tilde{Z},\tilde{T}\} = \partial_p\tilde{Z}\partial_q \tilde{T} - \partial_q\tilde{Z}\partial_p \tilde{T} = (\alpha_1\beta_2 -\alpha_2\beta_1)e^q,$$
$$P^2(\tilde{Z},\tilde{T}) = \Lambda^{12}\Lambda^{12}\partial_{pp}\tilde{Z} \partial_{qq}\tilde{T} +  \Lambda^{12}\Lambda^{21}\partial_{pq}\tilde{Z} \partial_{qp}\tilde{T} + 
 \Lambda^{21}\Lambda^{12}\partial_{qp}\tilde{Z} \partial_{pq}\tilde{T} +  $$ $$ +\Lambda^{21}\Lambda^{21}\partial_{qq}\tilde{Z} \partial_{pp}\tilde{T} = 0.$$
By analogy we have $$P^k(\tilde{Z},\tilde{T}) = 0, \forall k \geq 2.$$ Thus,
$$ i\tilde{Z} \star i\tilde{T} - i\tilde{T} \star i\tilde{Z} = \frac{1}{2i}[P^1(i\tilde{Z} , i\tilde{T}) - P^1(i\tilde{T}, i\tilde{Z})] = i(\alpha_1\beta_2 -\alpha_2\beta_1)e^q,$$ on one hand.

On the other hand, because $$[Z,T] = ZT-TZ = (\alpha_1\beta_2 - \alpha_2\beta_1) Y,$$ we have $$i\widetilde{[Z,T]} = i(\alpha_1\beta_2 -\alpha_2\beta_1)e^q =  i\tilde{Z} \star i\tilde{T} - i\tilde{T} \star i\tilde{Z} .$$ The proposition is hence proved.

Consequently, to each adapted chart $\psi$ in the sense of \cite{arnalcortet2}, we associate a $G$-covariant $\star$-product.

\begin{proposition}[see \cite{gutt}]
Let $\star$ be a formal differentiable $\star$-product on $C^\infty(M, \mathbb R)$, which is covariant under $G$. Then there exists a representation $\tau$ of $G$ in $\Aut N[[\nu]]$ such that 
$$\tau(g)(u \star v) = \tau(g)u \star \tau(g)v.$$
\end{proposition}

Let us denote by $\mathcal F_pu$ the partial Fourier transform \cite{meisevogt} of the function $u$ from the variable $p$ to the variable $x$, i.e.
$$\mathcal F_p(u)(x,q) := \frac{1}{\sqrt{2\pi}}\int_{\mathbb R} e^{-ipx} u(p,q)dp.$$ Let us denote by $ \mathcal F^{-1}_p(u) (x,q)$ the inverse Fourier transform. 
\begin{lemma}\label{lem3.1}
1. $\partial_p \mathcal F^{-1}_p(p.u) = i \mathcal F^{-1}_p(x.u)  $ ,

2. $ \mathcal F_p(v) = i \partial_x\mathcal F_p(v)  $ ,

3. $P^k(\tilde{Z},\mathcal F^{-1}_p(u)) = (-1)^k \beta e^q \frac{\partial^k\mathcal F^{-1}_p(u)}{\partial^kp}, \mbox{ with } k \geq 2.$
\end{lemma}

{\it Proof}. The first two formulas are well-known from theory of Fourier transforms. We reproduces them to locate notation.

1.  $\partial_p \mathcal F^{-1}_p(u) = \partial_p(\frac{1}{\sqrt{2\pi}}\int_{\mathbb R} e^{ipx} u(x,q)dx)= \frac{1}{\sqrt{2\pi}}\int_{\mathbb R} ix e^{ipx} u(x,q)dx =i \mathcal F^{-1}_p(x.u) $.

2.  $ i\partial_x \mathcal F_p(v) =i\partial_x(\frac{1}{\sqrt{2\pi}}\int_{\mathbb R} e^{-ipx} v(p,q)dp = i \frac{1}{\sqrt{2\pi}}\int_{\mathbb R} -ip e^{-ipx} v(p,q)dp = \frac{1}{\sqrt{2\pi}}\int_{\mathbb R} e^{-ipx} p v(p,q)dp =i \partial_x\mathcal F_p(p.v)  $ 

3. Remark that $\Lambda = \left(\begin{array}{cc} 0 & -1\\ 1 & 0 \end{array}\right)$ in the standard symplectic Darboux coordinates $(p,q)$ on the orbit $\Omega_+$ and we have had $\tilde{Z} = \alpha p + \beta e^q$, then 
$$P^2(\tilde{Z},\mathcal F^{-1}_p(u)) = \Lambda^{12}\Lambda^{12}\partial_{pp}\tilde{Z} \partial_{qq}\mathcal F^{-1}_p(u))  +  \Lambda^{12}\Lambda^{21}\partial_{pq}\tilde{Z} \partial_{qp}\mathcal F^{-1}_p(u))  + $$ $$
 \Lambda^{21}\Lambda^{12}\partial_{qp}\tilde{Z} \partial_{pq}\mathcal F^{-1}_p(u))   +\Lambda^{21}\Lambda^{21}\partial_{qq}\tilde{Z} \partial_{pp}\mathcal F^{-1}_p(u))  = (-1)^2\beta e^q \partial^2_{pp}\mathcal F^{-1}_p(u)).$$
By analogy we have $$P^k(\tilde{Z}, \mathcal F^{-1}_p(u) ) = (-1)^k\beta e^q \partial^k_{p\dots p}\mathcal F^{-1}_p(u)), \forall k \geq 3.$$ The lemma is therefore proved.

For each $Z \in \aff(\mathbb R)$, the corresponding Hamiltonian function is $\tilde{Z} =  \alpha p + \beta e^q $ and we can consider the operator $\ell_Z$ acting on dense subspace $L^2(\mathbb R^2, \frac{dpdq}{2\pi})^\infty$ of smooth functions by left $\star$-multiplication by $i \tilde{Z}$, i.e. $\ell_Z(u) = i\tilde{Z} \star u$. It is then continuated to the whole space $L^2(\mathbb R^2, \frac{dpdq}{2\pi})$. It is easy to see that, because of the relation in Proposition (\ref{3.1}), the correspondence $Z \in \aff(\mathbb R) \mapsto \ell_Z = i\tilde{Z} \star .$ is a representation of the Lie algebra $\aff(\mathbb R)$ on the space $N[[\frac{i}{2}]]$ of formal power series in the parameter $\nu = \frac{i}{2}$ with coefficients in $N = C^\infty(M,\mathbb R)$, see e.g. \cite{gutt} for more detail.

We study now the convergence of the formal power series. In order to do this, we look at the $\star$-product of $i\tilde{Z}$ as the $\star$-product of symbols and define the differential operators corresponding to $i\tilde{Z}$. It is easy to see that the resulting correspondence is a representation of $\mathfrak g $ by pseudo-differential operators. 

\begin{proposition}
For each $Z \in \aff(\mathbb R)$ and for each compactly supported $C^\infty$ function $u \in C^\infty_0(\mathbb R^2)$, we have 
$$\hat{\ell}_Z(u) := \mathcal F_p \circ \ell_Z \circ \mathcal F^{-1}_p(u) = \alpha (\frac{1}{2}\partial_q - \partial_x)u + i\beta e^{q -\frac{x}{2}}u.$$
\end{proposition}
{\it Proof}. 
For each $Z \in \mathfrak g = \aff(\mathbb R)$, we have
$$\hat{\ell}_Z(u) := \mathcal F_p \circ \ell_Z \circ \mathcal F^{-1}_p(u)  = \mathcal F_p(i\tilde{Z} \star \mathcal F^{-1}_p(u)) = i\mathcal F_p(\sum_{r \geq 0} \left(\frac{1}{2i}\right)^r P^r(\tilde{Z},\mathcal F^{-1}_p(u)).$$ Remark that $$P^1(\tilde{Z},\mathcal F^{-1}_p(u)) = \{\tilde{Z},\mathcal F^{-1}_p(u) \} = \alpha \partial_q \mathcal F^{-1}_p(u) - \beta e^q \partial_p \mathcal F^{-1}_p(u)$$ and applying Lemma (\ref{lem3.1}), we obtain:
$$i\mathcal F_p(\sum_{r\geq 0}\frac{1}{r!}\left(\frac{1}{2i}\right)^r P^r(\tilde{Z},\mathcal F^{-1}_p(u)) = $$
$$\begin{array}{cl}  &= i\mathcal F_p[(\alpha p + \beta e^q)\mathcal F_p^{-1}(u) + \frac{1}{2i}\alpha \partial_q \mathcal F_p^{-1}(u) - \frac{1}{2i}\beta e^q \partial_p\mathcal F^{-1}_p(u) +\\ 
 &+ \frac{1}{2!}\left(\frac{-1}{2i}\right)^2\beta e^q \partial_p^2\mathcal F^{-1}_p(u) + \dots + \frac{1}{n!}\left(\frac{-1}{2i}\right)^n\beta e^q \partial_p^n\mathcal F^{-1}_p(u) + \dots ] \\
 &= i[\alpha i \partial_x u + \beta e^qu + \frac{1}{2i}\alpha\partial_qu - \frac{1}{2i}\beta e^q \mathcal F_p(i\mathcal F^{-1}_p(x.u)) +\\
&+ \frac{1}{2!} \left(-\frac{1}{2i}\right)^2\beta e^q \mathcal F_p(i^2\mathcal F^{-1}_p(x^2.u)) + \dots +\frac{1}{n!} \left(-\frac{1}{2i}\right)^n\beta e^q \mathcal F_p(i^n\mathcal F^{-1}_p(x^n.u)) + \dots ]\\
 &= i[i\alpha \partial_x u + \frac{1}{2i}\alpha \partial_qu + \beta e^q u - \beta e^q \frac{x}{2} u +\\
&+ \frac{1}{2!}\beta e^q\left(\frac{x}{2}\right)^2 u + \dots + \frac{1}{n!}(-1)^n\beta e^q\left(\frac{x}{2}\right)^n u + \dots ]\\
 &= \alpha(\frac{1}{2}\partial_q - \partial_x)u + i\beta e^q[1 - \frac{x}{2}+ \frac{1}{2!}\left(\frac{x}{2}\right)^2 + \dots + (-1)^n \frac{1}{n!}\left(\frac{x}{2}\right)^n + \dots]\\
 &=  \alpha(\frac{1}{2}\partial_q - \partial_x)u + i\beta e^{q-\frac{x}{2}}u .
\end{array}$$
The proposition is therefore proved.

\begin{remark}{\rm
Setting new variables $s = q - \frac{x}{2}$, $t = q + \frac{x}{2}$, we have 
\begin{equation}
\hat{\ell}_Z(u) = \alpha\frac{\partial u}{\partial s} + i\beta e^s u,
\end{equation}
e.i. $$\hat{\ell}_Z = \alpha\frac{\partial }{\partial s} + i\beta e^s ,$$ which provides a representation of the Lie algebra $\aff (\mathbb R)$. 
}\end{remark}

\section{The associate irreducible unitary representations}

Our aim in this section is to exponentiate the obtained representation $\hat{\ell}_Z$ of the Lie algebra $\aff(\mathbb R)$ to the corresponding representation of the Lie group $\Aff_0(\mathbb R)$. We shall prove that the result is exactly the irreducible unitary representation $T_{\Omega_+}$ obtained from the orbit method or Mackey small subgroup method applied to this group $\Aff(\mathbb R)$.
Let us recall first the well-known list of all the irreducible unitary representations of the group of affine transformation of the real straight line.
\begin{theorem} [\cite{gelfandnaimark}]\label{4.1}
Every irreducible unitary representation of the group $\Aff(\mathbb R)$ of all the affine transformations of the real straight line, up to unitary equivalence, is equivalent to one of the pairwise nonequivalent representations:
\begin{itemize} 
\item the infinite dimensional representation $S$, realized in the space $L^2(\mathbb R^*, \frac{dy}{\vert y\vert})$, where $\mathbb R^* = \mathbb R \setminus \{0\}$ and is defined by the formula
$$(S(g)f)(y) := e^{iby}f(ay), \mbox{ where } g = \left(\begin{array}{cc} a & b\\ 0 & 1 \end{array}\right),$$
\item the representation $U^\varepsilon_\lambda$, where $\varepsilon = 0,1$, $\lambda \in \mathbb R$, realized in the 1-dimensional Hilbert space $\mathbb C^1$ and is given by the formula
$$U^\varepsilon_\lambda(g) = \vert a \vert^{i\lambda}(\sgn a)^\varepsilon .$$
\end{itemize}
\end{theorem}
Let us consider now the connected component $G= \Aff_0(\mathbb R)$. The irreducible unitary representations can be obtained easily from the orbit method machinery.
\begin{theorem}
The representation $\exp(\hat{\ell}_Z)$ of the group $G=\Aff_0(\mathbb R)$ is exactly the irreducible unitary representation $T_{\Omega_+}$ of $G=\Aff_0(\mathbb R)$ associated following the orbit method construction, to the orbit $\Omega_+$, which is the upper half-plane $\mathbb H \cong \mathbb R \rtimes \mathbb R^*$, i. e.
+$$(\exp(\hat{\ell}_Z)f)(y) = (T_{\Omega_+}(g)f)(y) = e^{iby}f(ay),\forall f\in L^2(\mathbb R^*, \frac{dy}{\vert y\vert}), $$ where  $g = \exp Z = \left(\begin{array}{cc} a & b \\ 0 & 1 \end{array}\right).$
\end{theorem}
{\it Proof}. Following the orbit method construction \cite{dndiep}, \cite{kirillov1}.
We choose an admissible Lie sub-algebra ${\mathfrak h} = \langle X \rangle$. Let us denote by $H$ the corresponding analytic subgroup of $G$ with Lie algebra $\mathfrak h$. The corresponding representation $\Ind_H^G\chi_F = \Ind_H^G\chi_{Y^*}$. The homogeneous space $H\setminus G$ is homeomorphic to $\mathbb R^* =\mathbb R \setminus \{0\}$ with the quasi-invariant measure $\frac{dy}{\vert y\vert}$. The corresponding representation $T_{\Omega_+}$ is given exactly by the same formula as the representation $S$ in the Theorem (\ref{4.1}). More precisely, for the element 
$$Z =\left(\begin{array}{cc} \alpha & \beta \\ 0 & 0 \end{array}\right) \in \mathfrak g = \aff(\mathbb R),$$
$$\exp Z = \exp \left(\begin{array}{cc} \alpha & \beta \\ 0 & 0 \end{array}\right)  = \left(\begin{array}{cc} a & b \\ 0 & 0 \end{array}\right) = \left\{\begin{array}{ll} \left(\begin{array}{cc} e^\alpha & \frac{\beta}{\alpha}(e^\alpha -1) \\ 0 & 0 \end{array}\right) &\mbox{if } \alpha \ne 0\\ \left(\begin{array}{cc} 1 & \beta \\ 0 & 1 \end{array}\right) &\mbox{if } \alpha = 0 \end{array}\right.$$ 
It is reasonable to simplify the notation, to consider the second case  
Remark that $y = e^q$ is the natural but non-canonical coordinate in $\mathbb R^* \cong H \setminus G$ we can write the induced representation obtained from the orbit method construction as
\begin{equation}T_{\Omega_+}(\exp Z)f(e^s) = \exp(i\frac{\beta}{\alpha}(e^{\alpha}-1)e^s)f(e^{\alpha + s}).\end{equation}
Therefore for the one-parameter subgroup $\exp(tZ), t\in \mathbb R$, we have the action formula
$$T_{\Omega_+}(\exp tZ)f(e^s) = \exp(i\frac{\beta}{\alpha}(e^{t\alpha}-1)e^s)f(e^{t\alpha + s}).$$
By a direct computation, we obtain
\begin{equation}\label{10} \frac{\partial}{\partial t}T_{\Omega_+}(\exp tZ)f(e^s) = \end{equation} 
$$\begin{array}{rl} &= i\frac{\beta}{\alpha}e^s\alpha e^{t\alpha} \exp(i\frac{\beta}{\alpha}(e^{t\alpha}-1)e^s)f(e^{ts+s}) + \exp(i\frac{\beta}{\alpha}(e^{t\alpha}-1)e^s) \alpha e^{t\alpha +s} \frac{\partial f}{\partial s}\\
 &=e^{t\alpha + s} \exp(i\frac{\beta}{\alpha}(e^{t\alpha}-1)e^s)[i\beta f(e^{t\alpha + s}) + \alpha \frac{\partial f}{\partial s}], \end{array}$$ on one hand.

On the other hand, we have
\begin{equation} \label{11} \hat{\ell}_Z T_{\Omega_+}(\exp tz)f(e^s) =  \end{equation}
$$\begin{array}{ll} &= (i\beta e^s +\alpha \frac{\partial}{\partial s})[\exp(i\frac{\beta}{\alpha}(e^{t\alpha}-1)e^s)f(e^{t\alpha + s})] \\
&= i\beta e^s \exp(i\frac{\beta}{\alpha}(e^{t\alpha}-1)e^s)f(e^{t\alpha + s}) +\\
&+ \alpha[i\frac{\beta}{\alpha}(e^{t\alpha} -1)e^s \exp(i\frac{\beta}{\alpha}(e^{t\alpha}-1)e^s)f(e^{t\alpha+s}) + \exp(i\frac{\beta}{\alpha}(e^{t\alpha}-1)e^s)e^{t\alpha+s}\frac{\partial f}{\partial s}]\\
&= e^{t\alpha+s}\exp(i\frac{\beta}{\alpha}(e^{t\alpha}-1)e^s)[i\beta f(e^{t\alpha+s}) + \alpha \frac{\partial f}{\partial s}].\end{array}$$
From (\ref{10}) and (\ref{11}) implies that 
$$\frac{\partial}{\partial t}T_{\Omega_+}(\exp(tZ))f(y) = \hat{\ell}_ZT_{\Omega_+}(\exp(tZ))f(y) .$$ Obviously, $$T_{\Omega_+}(\exp(tZ))f(y)|_{t=0} = f(y).$$ This means that $T_{\Omega_+}(\exp(tZ))f(y) $ is the unique solution of the Cauchy problem
$$\left\{ \begin{array}{rl} \frac{\partial}{\partial t}U(t,y) &= \hat{\ell}_Z U(t,y)\\ 
      U(0,y) &= \Id \end{array} \right.$$ 
This means also that $$\exp(\hat{\ell}_Z)f(y) \equiv T_{\Omega_+}(\exp Z)f(y).$$ The proof of the theorem is therefore achieved.

By analogy, we have also
\begin{theorem}
The representation $\exp(\hat{\ell}_Z)$ of the group $G=\Aff_0(\mathbb R)$ is exactly the irreducible unitary representation $T_{\Omega_-}$ of $G=\Aff_0(\mathbb R)$ associated following the orbit method construction, to the orbit $\Omega_-$, which is the lower half-plane $\mathbb H \cong \mathbb R \rtimes \mathbb R^*$, i. e.
$$(\exp(\hat{\ell}_Z)f)(y) = (T_{\Omega_-}(g)f)(y) = e^{iby}f(ay),\forall f\in L^2(\mathbb R^*, \frac{dy}{\vert y\vert}), $$ where  $g = \exp Z = \left(\begin{array}{cc} a & b \\ 0 & 1 \end{array}\right).$
\end{theorem}
\begin{remark}{\rm
1. We have demonstrated how all the irreducible unitary representation of the connected group of affine transformations could be obtained from deformation quantization. It is reasonable to refer to the algebras of functions on co-adjoint orbits with this $\star$-product as {\it quantum ones}.

2. In a forthcoming work, we shall do the same calculation for the group of affine transformations of the complex straight line $\mathbb C$. This achieves the description of {\it quantum $\overline{MD}$ co-adjoint orbits}, see \cite{dndiep} for definition of $\overline{MD}$ Lie algebras.
}\end{remark}

\end{document}